\documentclass[12pt]{amsart}
\usepackage{latexsym,amsmath, amscd, amssymb, amsthm, euscript,
mathrsfs}
\usepackage[all]{xy}
\usepackage{hyperref}

\newtheorem{thm}[subsection]{Theorem}

\newtheorem{cor}[subsection]{Corollary}
\newtheorem{lem}[subsection]{Lemma}
\newtheorem{lemma}[subsection]{Lemma}

\newtheorem{prop}[subsection]{Proposition}

\theoremstyle{definition}

\newtheorem{defn}[subsection]{Definition}

\theoremstyle{definition}

\theoremstyle{definition}

\newtheorem{rem}[subsection]{Remark}

\newtheorem{example}[subsection]{Example}
\numberwithin{equation}{subsection}

\newcommand{\Q}{\mathbb{Q}}

\newcommand{\V}{\mathbb{V}}

\newcommand{\mc}{\mathcal }

\newcommand{\Z}{\mathbb{Z}}

\newcommand{\Sp}{\text{\rm Spec}}
\newcommand{\mymargin}[1]{}

  \DeclareMathSymbol{A}{\mathalpha}{operators}{`A}%
  \DeclareMathSymbol{B}{\mathalpha}{operators}{`B}%
  \DeclareMathSymbol{C}{\mathalpha}{operators}{`C}%
  \DeclareMathSymbol{D}{\mathalpha}{operators}{`D}%
  \DeclareMathSymbol{E}{\mathalpha}{operators}{`E}%
  \DeclareMathSymbol{F}{\mathalpha}{operators}{`F}%
  \DeclareMathSymbol{G}{\mathalpha}{operators}{`G}%
  \DeclareMathSymbol{H}{\mathalpha}{operators}{`H}%
  \DeclareMathSymbol{I}{\mathalpha}{operators}{`I}%
  \DeclareMathSymbol{J}{\mathalpha}{operators}{`J}%
  \DeclareMathSymbol{K}{\mathalpha}{operators}{`K}%
  \DeclareMathSymbol{L}{\mathalpha}{operators}{`L}%
  \DeclareMathSymbol{M}{\mathalpha}{operators}{`M}%
  \DeclareMathSymbol{N}{\mathalpha}{operators}{`N}%
  \DeclareMathSymbol{O}{\mathalpha}{operators}{`O}%
  \DeclareMathSymbol{P}{\mathalpha}{operators}{`P}%
  \DeclareMathSymbol{Q}{\mathalpha}{operators}{`Q}%
  \DeclareMathSymbol{R}{\mathalpha}{operators}{`R}%
  \DeclareMathSymbol{S}{\mathalpha}{operators}{`S}%
  \DeclareMathSymbol{T}{\mathalpha}{operators}{`T}%
  \DeclareMathSymbol{U}{\mathalpha}{operators}{`U}%
  \DeclareMathSymbol{V}{\mathalpha}{operators}{`V}%
  \DeclareMathSymbol{W}{\mathalpha}{operators}{`W}%
  \DeclareMathSymbol{X}{\mathalpha}{operators}{`X}%
  \DeclareMathSymbol{Y}{\mathalpha}{operators}{`Y}%
  \DeclareMathSymbol{Z}{\mathalpha}{operators}{`Z}%

\def\CyrillicGuillemets{\DeclareFontEncoding{OT2}{}{}%
     \DeclareFontSubstitution{OT2}{wncyr}{m}{n}%
     \DeclareTextCommand{\guillemotleft}{OT1}{%
        {\fontencoding{OT2}\fontfamily{wncyr}\selectfont\char60}}%
     \DeclareTextCommand{\guillemotright}{OT1}{%
        {\fontencoding{OT2}\fontfamily{wncyr}\selectfont\char62}}}

\def\LasyGuillemets{%
     \DeclareTextCommand{\guillemotleft}{OT1}{\hbox{%
        \fontencoding{U}\fontfamily{lasy}\selectfont(\kern-0.20em(}}%
     \DeclareTextCommand{\guillemotright}{OT1}{\hbox{%
        \fontencoding{U}\fontfamily{lasy}\selectfont)\kern-0.20em)}}}
  \IfFileExists{ot2wncyr.fd}{\CyrillicGuillemets}{\LasyGuillemets}
   \DeclareTextSymbolDefault{\guillemotleft}{OT1}
   \DeclareTextSymbolDefault{\guillemotright}{OT1}
   \def\guill@spacing{\penalty\@M\hskip.8\fontdimen2\font
                               plus.3\fontdimen3\font
                               minus.8\fontdimen4\font}

\newcommand{\ra}{\rightarrow}

\newcommand{\LE}{{\textup{Lisse-{\'e}t}}}

\newcommand{\et}{\textup{{\'e}t}}

\newcommand{\X}{{\mc X}}
\newcommand{\NN}{{\mathbf N}}

\newcommand{\D}{{\mc D}}

\newcommand{\DD}{{\mathbf D}}

\DeclareMathOperator{\Hp}{{{}^p\mc H}}

\DeclareMathOperator{\Dp}{{{}^p\D}}
\DeclareMathOperator{\DDp}{{{}^p\DD}}

\DeclareMathOperator{\ext}{\mathcal{E}\it{xt}}

\DeclareMathOperator{\h}{{{}^p\mc H}}
\begin{document}

\title{Perverse sheaves on Artin stacks}

\author{Yves Laszlo and Martin Olsson}
\address{\'Ecole Polytechnique CMLS UMR 7640 CNRS F-91128 Palaiseau Cedex France}
 \email{laszlo@math.polytechnique.fr}
 \address{University of Texas at Austin
Department of Mathematics 1 University Station C1200 Austin, TX
78712-0257, USA}\email{molsson@math.utexas.edu}

\begin{abstract}In this paper we develop the theory of perverse
sheaves on Artin stacks continuing  the study in
~\cite{Las-Ols05-1} and~\cite{Las-Ols05-2}.\end{abstract}

\maketitle

\section{Introduction}

In this third paper in our series on Grothendieck's six operations for {\'e}tale sheaves on stacks, we define the perverse $t$--structure on the derived category of {\'e}tale sheaves (with either finite or adic coefficients). This generalizes the $t$--structure defined in the case of schemes in \cite{BBD82} (but note also that in this paper we consider unbounded schemes which is not covered in loc. cit.).

By \cite[3.2.4]{BBD82} the perverse sheaves on a scheme form a stack with respect to the smooth topology.  This enables one to define the notion of a perverse sheaf on any Artin stack (using also the unbounded version of the gluing lemma \cite[3.2.4]{BBD82} proven in \cite[2.3.3]{Las-Ols05-1}).   The main contribution of this paper is to define a $t$--structure on the derived category whose heart is this category of perverse sheaves.  In fact it is shown in \cite[4.2.5]{BBD82} that pullback along smooth morphisms of schemes is an exact functor (up to a shift) with respect to the perverse $t$--structure.  This implies that the definition of the functor $\h^0$ ($\tau _{\geq 0}\tau _{\leq 0}$ with respect to the perverse $t$--structure) for stacks is forced upon us from the case of schemes and the gluing lemma.  We verify in this paper that the resulting definitions of the subcategories $\Dp^{\geq j}$ and $\Dp^{\leq j}$ of the derived category of {\'e}tale sheaves (in either finite coefficient or adic coefficient case) define a $t$--structure. 

\begin{rem} The reader should note that unlike the case of schemes (Beilinson's theorem \cite{Bei87}) the derived category of the abelian category of perverse sheaves is not equivalent to the derived category of sheaves on the stack.  An explicit example suggested by D. Ben-Zvi is the following: Let $\mc X = B\mathbb{G}_m$ over an algebraically closed field $k$.  The category of perverse $\mathbb{Q}_l$--sheaves on $\mc X$ is the equivalent to the category of $\mathbb{G}_m$--equivariant perverse sheaves on $\Sp (k)$ as defined in \cite[III.15]{KW}.  In particular, this is a semisimple category.  In particular, if $\mc D$ denotes the derived category of perverse sheaves we see that for two objects $V$ and $W$ the groups $\text{Ext}^{i}_{\mc D}(V, W)$ are zero for all $i>0$.  On the other hand, we have $H^{2i}(\mc X, \mathbb{Q}_l) \neq 0$ for all $i\geq 0$.
\end{rem}

\begin{rem} The techniques used in this paper can also be used to define perverse sheaves of $\mc D$--modules on complex analytic stacks.
\end{rem}

Throughout the paper we work over a ground field $k$ and write $S = \Sp (k)$.

\section{Gluing of $t$--structures}

For the convenience of the reader, we review the key result \cite[1.4.10]{BBD82}.

Let $\mc D$, $\mc D_U$, and $\mc D_F$ be three triangulated categories with exact functors
$$
\begin{CD}
\mc D_F@>i_*>> \mc D@>j^*>> \mc D_U.
\end{CD}
$$
Write also $i_!:= i_*$ and $j^!:= j^*$.  Assume the following hold:
\begin{enumerate}\label{gluingcond}
\item [(i)] The functor $i_*$ has a left adjoint $i^*$ and a right adjoint $i^!$.
\item [(ii)] The functor $j^*$ has a left adjoint $j_!$ and a right adjoint $j_*$.
\item [(iii)] We have $i^!j_* = 0$.
\item [(iv)] For every object $K\in \mc D$ there exists  a morphism $d:i_*i^*K\rightarrow j_!j^*K[1]$ (resp. $d:j_*j^*K\rightarrow i_*i^!K[1]$) such that the induced triangle
$$
j_!j^*K\rightarrow K\rightarrow i_*i^*K\rightarrow j_!j^*K[1]
$$
$$
(\text{resp. } i_*i^!K\rightarrow K\rightarrow j_*j^*K\rightarrow i_*i^!K[1])
$$
is distinguished.
\item [(v)] The adjunction morphisms $i^*i_*\rightarrow \text{id}\rightarrow i^!i_*$ and $j^*j_*\rightarrow \text{id}\rightarrow j^*j_!$ are all isomorphisms.
\end{enumerate}

The main example we will consider is the following:

\begin{example}\label{mainexample}
Let $\X$ be an algebraic stack locally of finite type over a field $k$, and let $i:\mc F\hookrightarrow \X$ be a closed substack with complement $j:\mc U\hookrightarrow \X$.  Let $\Lambda $ be a complete discrete valuation ring of residue characteristic prime to $\text{char}(k)$, and for an integer $n$ let $\Lambda _n$ denote $\Lambda /\mathfrak{m}^{n+1}$.

Fix an integer $n$, and let 
$\mc D$ (resp. $\mc D_{\mc U}$, $\mc D_{\mc F}$) denote the bounded derived category $\D_c^b(\X, \Lambda _n)$ (resp. $\D_c^b(\mc U, \Lambda _n)$, $\D_c^b(\mc F, \Lambda _n)$), and let $i_*:\D_{\mc F}\rightarrow \D$ and $j^*:\D\rightarrow \D_{\mc U}$ be the usual pushforward and pullback functors.  By the theory developed in \cite{Las-Ols05-1} conditions (i)-(v) hold.  

We can also consider adic sheaves.  Let  $\mc D$ (resp. $\mc D_{\mc U}$, $\mc D_{\mc F}$) denote the bounded derived category $\DD_c^b(\X, \Lambda )$ (resp. $\DD_c^b(\mc U, \Lambda )$, $\DD_c^b(\mc F, \Lambda )$) of $\Lambda $--modules on $\X$ (resp. $\mc U$, $\mc F$) constructed in \cite{Las-Ols05-2}. We then again have functors
$$
\begin{CD}
\D_{\mc F}@>i_*>> \D@>j^*>> \D_{\mc U}.
\end{CD}
$$
Conditions (i)-(iii) hold by the results of \cite{Las-Ols05-2}, and condition (v) holds by base change to a smooth cover of $\X$ and the case of schemes.

To construct the distinguished triangles in (iv) recall that $\DD_c^b(\X, \Lambda )$ is constructed as a quotient of the category $\D_c^b(\X^{\mathbb{N}}, \Lambda _\bullet )$ (the derived category of projective systems of $\Lambda _n$--modules), and similarly for $\DD_c^b(\mc U, \Lambda )$ and $\DD_c^b(\mc F, \Lambda )$.  All the functors in (i)-(v) are then obtained from functors defined already on the level of the categories $\D_c^b(\X^{\mathbb{N}}, \Lambda _\bullet )$, $\D_c^b(\mc U^{\mathbb{N}}, \Lambda _\bullet )$, and $\D_c^b(\mc F^{\mathbb{N}}, \Lambda _\bullet )$.  In this case the first distinguished triangle in (iv) is constructed by the same reasoning as in \cite[4.9]{Las-Ols05-1} for the finite case, and the second distinguished triangle is obtained by duality.
\end{example}

Returning to the general setup of the beginning of this section,
suppose given $t$--structures $(\mc D_F^{\leq 0}, \mc D_F^{\geq 0})$ and $(\mc D_U^{\leq 0}, \mc D_U^{\geq 0})$ on $\mc D_F$ and $\mc D_U$ respectively and define
$$
\mc D^{\leq 0}:= \{K\in \mc D|j^*K\in \mc D^{\leq 0}_U \text{ and } i^*K\in \mc D_F^{\leq 0}\}
$$
$$
\mc D^{\geq 0}:= \{K\in \mc D|j^*K\in \mc D^{\geq 0}_U \text{ and } i^!K\in \mc D_F^{\geq 0}\}.
$$

\begin{thm}[{\cite[1.4.10]{BBD82}}]\label{gluetstructure} The pair $(\mc D^{\leq 0}, \mc D^{\geq 0})$ defines a $t$--structure on $\mc D$.
\end{thm}

\section{Review of the perverse $t$--structure for schemes}

Let $k$ be a field and $X/k$ a scheme of finite type.  Let $\Lambda $ be a complete discrete valuation ring and for every $n$ let $\Lambda _n$ denote the quotient $\Lambda /\mathfrak{m}^{n+1}$ so that $\Lambda = \varprojlim \Lambda _n$.  Assume that the characteristic $l$ of $\Lambda _0$ is invertible in $k$.  

For every $n$, we can define the perverse $t$--structure $(\Dp^{\leq 0}(X, \Lambda _n), \Dp^{\geq 0}(X, \Lambda _n))$ on $\mc D^b_c(X, \Lambda _n)$ (in this paper we consider only the middle perversity) as follows:
\begin{enumerate}
\item [{}] A complex $K\in \mc D^b_c(X, \Lambda _n)$ is in $\Dp ^{\leq 0}(X, \Lambda _n)$ (resp. $\Dp ^{\geq 0}(X, \Lambda _n)$) if for every point $x\in X$ with inclusion $i_x:\Sp (k(x))\rightarrow X$  and $j>-\text{dim}(x)$ (resp. $j<-\text{dim}(x)$) we have $\mc H^j(i^*_xK) = 0$ (resp. $\mc H^j(i^!_xK) = 0$)\footnote{As usual, $i_x^!K,i_x^*K$ denotes $(i_{\bar x}^!K)_x, (i_{\bar x}^*K)_x$ where $i_{\bar x}$ is the closed immersion ${\bar x}_{red}\hookrightarrow X$.}.
\end{enumerate}

As explained in \cite[2.2.11]{BBD82} this defines a $t$--structure on $\D_c^b(X, \Lambda _n)$: The \emph{perverse} $t$--structure.

The same technique can be used in the adic case.  We explain this in more detail since it is not covered in  detail in the literature.    As before let  $\DD_c^b(X, \Lambda )$ denote the bounded derived category of $\Lambda $--modules constructed in \cite{Las-Ols05-2}.  Let $\text{Mod}_c(X, \Lambda )$ denote the heart of the standard $t$--structure on $\DD_c(X, \Lambda )$.  In the language of \cite[3.1]{Las-Ols05-2} the category $\text{Mod}_c(X, \Lambda )$ is the quotient of the category of $\lambda $--modules on $X$ by almost zero systems.  For every integer $j$ there is then a natural functor
$$
\mc H^j:\DD_c^b(X, \Lambda )\rightarrow \text{Mod}_c(X, \Lambda ).
$$

We then define categories $(\DDp^{\leq 0}(X, \Lambda ), \DDp^{\geq 0}(X, \Lambda ))$ by the following condition:
\begin{enumerate}
\item [{}] A complex $K\in \DD^b_c(X, \Lambda )$ is in $\DDp ^{\leq 0}(X, \Lambda )$ (resp. $\DDp ^{\geq 0}(X, \Lambda )$) if for every point $x\in X$ with inclusion $i_x:\Sp (k(x))\rightarrow X$  and $j>-\text{dim}(x)$ (resp. $j<-\text{dim}(x)$) we have $\mc H^j(i^*_xK) = 0$ (resp. $\mc H^j(i^!_xK) = 0$). 
\end{enumerate}

\begin{prop}\label{3.1b} This defines a $t$--structure on $\DD_c^b(X, \Lambda )$.
\end{prop}
\begin{proof} The only problem comes from perverse
truncation. Recall that an adic sheaf $(M_n)$ is \emph{smooth} if all $M_n$
are locally constant, or, what's amount to the same, if $M_1$ is locally
constant. We say that a complex $K\in \DD_c^b(X, \Lambda )$ is \emph{smooth} if $\mc H^j(K)$ is represented by a smooth adic sheaf and is zero for almost all $j$.

We say that $X$ is \emph{essentially smooth} if $(X\otimes _k\bar k)_{\text{red}}$ is smooth over $\bar k$.  If $X$ is essentially smooth of dimension $d$ and $K\in \DD^b_c(X, \Lambda )$ is smooth, then for any point $x\in X$ of codimension $s$ we have $i_x^!K = i_x^*K(s-d)[2(s-d)]$ (see for example \cite[p. 62]{SGA5}).

We now prove by induction on $\dim(X)$ that the third axiom for a $t$--structure holds.  Namely, that for any $K\in \DD^b_c(X, \Lambda )$ there exists a distinguished triangle
\begin{equation}\label{ptriangle}
{}^p\tau_{\leq 0}K\ra
K\ra{}^p\tau_{> 0}K
\end{equation}
with ${}^p\tau _{\leq 0}K\in \DDp^{\leq 0}(K, \Lambda )$ and ${}^p\tau _{>0}K\in \DDp^{>0}(K, \Lambda )$.

For $\dim(X)=0$, it is clearly true. 
For the inductive step let $d$ be the dimension of $X$ and assume the result holds for schemes of dimension $<d$.
Let  $K\in\DD^b_c(X,\Lambda)$ be a complex and
choose some essentially smooth dense open subset $U$ of $X$ on which
$K$ is smooth. Then, the class of $\tau_{\leq -\dim(U)}K_{|U}$ (resp.
$\tau_{> -\dim(U)}K_{|U}$) (truncation with respect to the usual $t$-structure on $U$) belongs to $\DDp
^{\leq 0}(U, \Lambda )$ (resp. $\DDp ^{> 0}(U, \Lambda )$) and therefore
the usual distinguished triangle
$$\tau_{\leq -\dim(U)}K_{K|U}\ra K_{|U}\ra\tau_{> -\dim(U)}K_{|U}$$
defines the required perverse distinguished triangle on $U$ by the formula
above. The complement $F=X-U$ has dimension $<\dim(X)$. By induction
hypothesis, the conditions above define a $t$-structure on $F$ and therefore
one gets a distinguished triangle $${}^p\tau_{\leq 0}K_{|F}\ra
K_{|F}\ra{}^p\tau_{> 0}K_{|F}$$
on $F$ . By \ref{gluetstructure} we can glue the trivial $t$--structure on $U$ and the perverse $t$--structure on $F$ to a $t$--structure on $\DD_c^b(X, \Lambda )$.  It follows that one can glue the distinguished triangles on $U$ and $F$ to a distinguished triangle \ref{ptriangle} which  gives the third axiom.\end{proof}

\begin{rem} One can also prove the proposition using stratifications as in \cite{BBD82}.
\end{rem}

The perverse $t$--structures on $\D_c^b(X, \Lambda _n)$ and $\DD_c^b(X, \Lambda )$ extend naturally to the unbounded derived categories $\D_c(X, \Lambda _n)$ and $\DD_c(\X, \Lambda )$. Let $\D$ denote either of these triangulated categories.  For $K\in \D$ and $a\leq b$ let $\tau _{[a,b]}K$ denote $\tau _{\geq a}\tau _{\leq b}K$.  The perverse $t$--structure defines a functor
$$
\h^0:\D^b\rightarrow \D^b.
$$

\begin{lem} There exists integer $a<b$ such that for any $K\in \D^b$ we have $\h^0(K) = \h^0(\tau _{[a,b]}K)$.
\end{lem}
\begin{proof} Consideration of the distinguished triangles
$$
\tau _{\leq b}K\rightarrow K\rightarrow \tau _{>b}K\rightarrow \tau _{\leq b}K[1]
$$
and
$$
\tau _{<a}K\rightarrow \tau _{\leq b}K\rightarrow \tau _{[a,b]}K\rightarrow \tau _{<a}K[1]
$$
implies that it suffices to show that there exists integers $a<b$ such that for $K$ in either $\D^{<a}$ or $\D^{>b}$ we have $\h^0(K) = 0$.  By the definition of perverse sheaf we can take $a$ to be any integer smaller than $-\text{dim}(X)$.

To find the integer $b$, note that since the dualizing sheaf for a scheme of finite type over $k$ has finite quasi--injective dimension \cite[I.1.5]{SGA5} and \cite[7.6]{Las-Ols05-2}.  It follows that there exists a constant $c$ such that for any integer $b$, point $x\in X$, and $K\in \D^{>b}$ we have $i_x^!K\in \D^{>b+c}$.  Thus we can take for $b$ any integer greater than $-\text{dim}(X)-c$.
\end{proof}

Choose integers $a<b$ as in the lemma, and define
$$
\h^0:\D\rightarrow \D^b, \ \ K\mapsto \h^0(\tau _{[a,b]}K).
$$
One sees immediately that this does not depend on the choice of $a<b$.  Define $\D^{\leq 0}$ (resp. $\D^{\geq 0}$) to be the full subcategory of $\D$ of complexes $K\in \D$ with $\h^j(K) := \h^0(K[j])=0$ for $j\leq 0$ (resp. $j\geq 0$).
The argument in \cite[2.2.1]{BBD82} (which in turn is based in \cite[2.1.4]{BBD82}) shows that this defines a $t$--structure on $\D$.

\section{The perverse $t$--structure for stacks of finite type}

Let $\X/k$ be an algebraic stack of finite type.  Let $\D(\X)$ denote either $\D_c(\X, \Lambda _n)$ or $\DD_c(\X, \Lambda )$.

Fix  a smooth surjection $\pi :X\rightarrow \X$ with $X$ of finite type, and define $\Dp^{\leq 0}(\X)$ (resp. $\Dp^{\geq 0}(\X)$) to be the full subcategory of objects $K\in \D(\X)$ such that $\pi ^*K[d_\pi ]$ is in $\Dp  ^{\leq 0}(X)$ (resp. $\Dp^{\geq 0}(X)$), where $d_\pi $ denotes the relative dimension of $X$ over $\X$ (a locally constant function on $X$).

\begin{lem} The subcategories $\Dp^{\leq 0}$ and $\Dp^{\geq 0}$ of $\D$ do not depend on the choice of $\pi :X\rightarrow \X$.
\end{lem}
\begin{proof} It suffices to show that if $f:Y\rightarrow X$ is a smooth surjective morphism of schemes of relative dimension $d_f$ (a locally constant function on $Y$), then $K\in \D(X)$ is in $\Dp^{\leq 0}(X)$ (resp. $\Dp^{\geq 0}(X)$) if and only if $f^*K[d_f]$ is in $\Dp ^{\leq 0}(Y)$ (resp. $\Dp ^{\geq 0}(Y)$).  For this note that by \cite[4.2.4]{BBD82} the functor $f^*[d_f]$ is exact for the perverse $t$--structures.  This implies that $K$ is in $\Dp ^{\leq 0}(X)$ (resp. $\Dp ^{\geq 0}(X)$) only if $f^*K[d_f]$ is in $\Dp ^{\leq 0}(Y)$ (resp. $\Dp ^{\geq 0}(Y)$).

For the other direction, note that if $f^*K[d_f]$ is in $\D^{\leq 0}(Y)$ (resp. $\D^{\geq 0}(Y)$) then for any integer $i>0$ (resp. $i<0$) we have
$$
f^*\h^i(K)[d_f] = \h^i(f^*K[d_f]) = 0.
$$
Since $f$ is surjective it follows that $\h^i(K) = 0$ for all $i>0$ (resp. $i<0$).
\end{proof}

\begin{thm} The subcategories $(\Dp^{\leq 0}(\X), \Dp ^{\geq 0}(\X))$ define a $t$--structure on $\D(\X)$.
\end{thm}
\begin{proof}
Exactly as in the proof of \ref{3.1b} using noetherian induction and gluing of $t$--structures \ref{gluetstructure} one shows that $(\Dp^{\leq 0}, \Dp^{\geq 0})$ define by restriction a $t$--structure on $\D^b(\X)$ (again the only problem is the third axiom for a $t$--structure since the other two can be verified locally).

The same argument used in the schematic case then extends this $t$--structure to the unbounded derived category $\D(\X)$.
\end{proof}

\section{The perverse $t$--structure for stacks locally of finite type}

Assume now that  $\mc X$ is a stack \emph{locally} of finite type
over $S$.  We consider either finite coefficients or the adic case and write just $\D(\X)$ for the corresponding derived categories $\D_c(\X, \Lambda _n)$ or $\DD_c(\X, \Lambda ).$

Define subcategories $(\Dp^{\leq 0}(\X), \Dp^{\geq 0}(\X))$ of $\D(\X)$ by the condition that $K\in \D(\X)$ is in $\Dp ^{\leq 0}(\X)$ (resp. $\Dp^{\geq 0}(\X)$) if and only if for every open substack $\mc U\subset \X$ of finite type over $k$ the restriction of $K$  to $\mc U$ is in $\Dp ^{\leq 0}(\mc U)$ (resp. $\Dp ^{\geq 0}(\mc U)$.

\begin{thm}\label{5.1thm} The subcategories $(\Dp^{\leq 0}(\X), \Dp^{\geq 0}(\X))$ define a $t$--structure on $\D(\X)$.
\end{thm}
\begin{proof}
The first two axioms for a $t$--structure follow immediately from the definition.  We now show the third axiom.
Write $\mc X$ as a filtering union of open
substacks $\mc
 X_i\subset \mc X$ of finite type.  Let $j_i:\mc X_i\hookrightarrow \mc
 X$ be the open immersion.  Then for any $M\in D_c(\mc X)$, we have for every $i$ a
 distinguished triangle

\begin{equation}\label{Ti}
     M_{i, \leq 0}\rightarrow M|_{\mc X_i}\rightarrow M_{i, \geq 1}
\end{equation} where $M_{i, \leq 0}\in\DDp^{\leq 0}(\X_i)$ and $M_{i, \leq 0}$ in
 $\DDp^{\geq 1}(\X_i)$.
By the uniqueness statement in \cite[1.3.3]{BBD82} this implies
 that the formation of this sequence is compatible with restriction to
 smaller $\mc X_i$.  Since $j_i^*=j_i^!$ for open immersions, we then get
 a sequence
 $$
 j_{i!}M_{i, \leq 0}\rightarrow j_{i+1, !}M_{i+1, \leq 0}\rightarrow
 \cdots.
 $$
 Define $M_{\leq 0}$ to be the homotopy colimit of this sequence.  There
 is a natural map $M_{\leq 0}\rightarrow M$ and take $M_{\geq 1}$ to be
 the cone. The following lemma implies that the third axiom holds and hence proves \ref{5.1thm}.
\end{proof}

 \begin{lemma}\label{triangle-pervers} For any $i$, the restriction of the distinguished triangle
 \begin{equation}\label{T}
M_{\leq 0}\rightarrow M\rightarrow M_{\geq 1}
\end{equation} $\X_i$ is isomorphic to \ref{Ti}. In particular, $M_{\leq 0}\in \Dp^{\leq 0}$ and $M_{\geq 1}\in \Dp^{\geq 1}$.
\end{lemma}

\begin{proof} Let $i_0$ be any nonnegative integer. By~\cite[1.7.1]{Nee01}, one has
 a distinguished triangle

$$\bigoplus_{i\geq i_0} j_{i!}M_{i, \leq 0}\rightarrow\bigoplus_{i\geq i_0} j_{i!}M_{i, \leq
0}\rightarrow M_{\leq 0}.$$

Because $j_i^*$ is exact and commutes with direct sums, one gets
by restriction a distinguished triangle

$$\bigoplus_{i\geq i_0} M_{i_0, \leq 0}\rightarrow\bigoplus_{i\geq i_0}M_{i_0, \leq
0}\rightarrow M_{\leq 0|\X_{i_0}}.$$ where the inductive system is
given by the identity morphism of $M_{i_0,\leq 0}$.
By\cite[1.6.6]{Nee01}, one gets $M_{\leq 0|\X_{i_0}}=M_{i_0,\leq
0}$. The lemma follows.\end{proof}

We define the \emph{perverse $t$--structure} on $\D$ to be the $t$--structure given by \ref{5.1thm}.
 By the very
 definition, it coincides with the usual one if $\X$ is a scheme.. A complex in the {heart} of the perverse
 t-structure is by definition a perverse sheaf.

\begin{rem} By \cite[1.3.6]{BBD82}, the category of perverse sheaves a stack $\X$ is an abelian category. 
\end{rem}

\begin{rem} If we work with $\Lambda _0$--coefficients, then it follows from the case of schemes that Verdier duality interchanges the categories $\Dp^{\leq }(\X, \Lambda _0)$ and $\Dp ^{\geq 0}(\X, \Lambda _0)$.  For other coefficients this does not  hold due to the presence of torsion.
\end{rem}

 \begin{rem}\label{perv-loc}If the normalized complex $P$ is perverse on $\X$ and
 $U\ra\X$ is in $\LE(\X)$, then $P_{U,n}\in\DD^b(U_\et,\Lambda_n)$ is perverse on $U_\et$. In
 particular, one has $\ext^i(P_{U,n},P_{U,n})=0$ if $i<0$. By the
 gluing lemma, perversity is a local condition for the lisse-{\'e}tale
 topology. For instance, it follows  that the
 category
 perverse sheaf on $\X=[X/G]$ ($X$ is a scheme of finite type acting on by an
 algebraic group $G$) is equivalent to the category of
 $G$-equivariant perverse sheaves on $X$\footnote{See~\cite[III.15]{KW}}.\end{rem}

In the case of finite coefficients, one can also define $\h^0$ by
gluing. Let us consider a diagram
\begin{equation}\label{UV}\xymatrix{V\ar[rr]^\sigma\ar[rd]_v&&U\ar[ld]^u\\
&\X}\end{equation} with a $2$-commutative triangle and
$u,v\in\LE(\X)$ of relative dimension $d_u,d_v$.  Let $R$ be a Gorenstein ring of dimension
$0$.

\begin{lemma}\label{pHOfini}
Let $K\in\D^b_c(\X,R)$. There exists a unique
$\h^0(K)\in\D^b(\X,R)$ such that
$$[\h^0(K)]_U[d_u]=\h^0(K_U[d_u])\in\D^b_c(U_\et,R)$$ (functorially).\end{lemma}

\begin{proof} Because $\h^0(K_U)$ is perverse, one has
by \cite[2.1.21]{BBD82}
$$\ext^i(\h^0(K_U[d_u]),\h^0(K_U[d_u]))=0\text{ for }i<0.$$ Let $W=U\times_\X V$ which is an algebraic space.

 Assume for
simplicity that $W$ is even a scheme (certainly of finite type
over $S$). One has a commutative diagram with cartesian square

$$\xymatrix{W\ar[r]^{\tilde v}\ar[d]_{\tilde u}&U\ar[d]^u\\V\ar@/^1pc/[u]^s
\ar[r]_v\ar[ru]^\sigma&\X}$$ with ${\tilde u},{\tilde v}$ smooth
of relative dimension $d_u,d_v$ and $s$ is a graph section. In
particular, ${\tilde u}^*[d_u]$ and ${\tilde v}^*[d_v]$ are
$t$-exact (for the perverse $t$-structure) by~\cite[4.2.4]{BBD82}.

Therefore, we get
\begin{eqnarray*}
  {\tilde v}^*\h^0(K_U[d_u])[-d_u] &=& {\tilde v}^*[d_v]\h^0(K_U[d_u])[-d_u-d_v] \\
   &=& \h^0({\tilde v}^*K_U[d_u+d_v])[-d_u-d_v] \\
   &=& \h^0(K_W[d_u+d_v])[-d_u-d_v] \\
   &=& {\tilde u}^*\h^0(K_V[d_v])[-d_v] \\
\end{eqnarray*}
Pulling back by $s$, we get
$$\h^0(K_V[d_V])[-d_v]=s^*{\tilde u}^*\h^0(K_V[d_V])[-d_v]=s^*{\tilde v}^*\h^0(K_U[d_u])[-d_u]=
\sigma^*\h^0(K_U[d_u])[-d_u].$$

The lemma follows from \cite[2.3.3]{Las-Ols05-1}.\end{proof}

\begin{rem} It follows from the construction of the perverse $t$--structure on $\D_c(\X, R)$ that the above defined functor $\h^0$ agrees with the one defined by the perverse $t$--structure.
\end{rem}

\section{Intermediate extension}

Let $\X$ be an algebraic $k$--stack of finite type, and let $i:\mc Y\hookrightarrow \X$ be a closed substack with complement $j:\mc U\hookrightarrow \X$.  For a perverse sheaf $P$ on $\mc U$ we define the \emph{intermediate extension}, denoted $j_{!*}P$,  to be the image in the abelian category of perverse sheaves on $\X$ of the morphism
$$
\h^0(j_!P)\rightarrow \h^0(j_*P).
$$

\begin{lem}\label{extlem} The perverse sheaf $j_{!*}P$ is the unique perverse sheaf with $j^*(j_{!*}P) = P$ and $\h^0(i^*(j_{!*}P)) = 0$.
\end{lem}
\begin{proof}
Let us first verify that $j_{!*}P$ has the indicated properties. Since $j$ is an open immersion, the functor $j^*$ is $t$--exact and hence the first property $j^*j_{!*}P = P$ is immediate. The equality $\h^0(i^*(j_{!*}P)) = 0$ follows from \cite[1.4.23]{BBD82}.

Let $F$ be a second perverse sheaf with these properties.  Then $j^*F=P$ defines a morphism $j_!P\rightarrow F$ which since $j_!$ is right exact for the perverse $t$--structure (this follows immediately from \cite[2.2.5]{BBD82} and a reduction to the case of schemes) factors through a morphism $\h^0(j_!P)\rightarrow F$.  Adjunction also defines a morphism $F\rightarrow j_*P$ which since $j_*$ is left exact for the perverse $t$--structure (again by loc. cit.) defines a morphism $F\rightarrow \h^0(j_*P)$.  It follows that the morphism $\h^0(j_!P)\rightarrow \h^0(j_*P)$ factors through $F$ whence we get a morphism $\rho :j_{!*}P\rightarrow F$ of perverse sheaves.  The kernel and cokernel of this morphism is a perverse sheaf supported on $\mc Y$.  The assumption $\h^0(i^*F) = \h^0(i^*j_{!*}P) = 0$ then implies that the kernel and cokernel are zero.
\end{proof}

\begin{lemma}\label{*!local} 
Let $f:X\rightarrow \X$ be a smooth morphism of relative dimension $d$ with $X$ a scheme.  Let 
$$
\begin{CD}
Y@>i'>> X@<j'<< U
\end{CD}
$$
be the pullbacks of $\mc Y$ and $\mc U$.  Then
 $f^*[d]j_{!*}=j'_{!*}f^*[d]$.\end{lemma}

\begin{proof} Let $P$ be a perverse sheaf on $\mc U$ and let $\bar P$ denote $j_{!*}P$.  The functor $f^*[d]$ is t-exact, and hence preserves perversity. It follows that  $\bar
P'=f^*[d]\bar P$ is perverse and is an extension of the perverse
sheaf  $P'=f^*[d]P$ (in particular the statement of the lemma
makes sense!). By the uniqueness in \ref{extlem} it suffices to show that
$\Hp^0(i'^*\bar P')=0.$ But, keeping in mind that $f^*[d]$ commutes with $\Hp^0$, the
first point is for instance a consequence of smooth base change. \end{proof}

\begin{rem} In the case of finite coefficients, one can also define the intermediate extension using \ref{*!local} and gluing.
\end{rem}

\section{Gluing perverse sheaves}

In this section we work either with finite coefficients or with adic coefficients.

Let $\mc X$ be a stack locally of finite type over $k$, and define a fibered category $\mc P$ (\emph{not} in groupoids) on $\LE(\X)$ by
$$
U\mapsto (\text{category of perverse sheaves on $U$}).
$$

\begin{prop}\label{perv-stack} The fibered category $\mc P$ is a stack and the natural functor
$$
(\text{\rm perverse sheaves on $\X$})\rightarrow \mc P(\X)
$$
is an equivalence of categories.
\end{prop}

\begin{proof} For a smooth surjective morphism of stacks $f:\mc Y\rightarrow \mc X$ let $\text{Des}_{\mc Y/\mc X}$ denote the category of pairs $(P, \sigma )$, where $P$ is a perverse sheaf on $\mc Y$ and $\sigma :\text{pr}_1^*P\rightarrow \text{pr}_2^*P$ is an isomorphism over $\mc Y\times _{\mc X}\mc Y$ satisfying the usual cocycle condition on $\mc Y\times _{\X}\mc Y\times _\X\mc Y$.  To prove the proposition it suffices to show that the natural functor
\begin{equation}\label{perversefunctor}
(\text{perverse sheaves on $\X$})\rightarrow \text{Des}_{\mc Y/\mc X}
\end{equation}
is an equivalence of categories.

Now if $P$ and $P'$ are perverse sheaves on a stack, then $\ext^i(P, P') = 0$ for all $i<0$. Indeed this can be verified locally where it follows from the first axiom of a $t$--structure. That \ref{perversefunctor} is an equivalence in the finite coefficients case then follows from the gluing lemma \cite[2.3.3 and 2.3.4]{Las-Ols05-1}.

For the adic case, note that by the discussion in \cite[\S 5]{Las-Ols05-2} if $P$ and $P'$ are two perverse sheaves on a stack $\X$ with normalized complexes $\hat P$ and $\hat P'$ then
$$
\text{Ext} ^i_{\D_c(\X^{\mathbb{N}}, \Lambda _\bullet )}(\hat P, \hat P') = \text{Ext}^i_{\DD_c(\X, \Lambda )}(P, P'),
$$
where $\X^{\mathbb{N}}$ denotes the topos of projective systems of sheaves on $\LE (\X)$ and $\Lambda _\bullet $ denotes $\varprojlim \Lambda _n$.  It follows that for any object $(P, \sigma )\in \text{Des}_{\mc Y/\mc X}$ we have $\ext^i_{\D_c(\mc Y^{\mathbb{N}}, \Lambda _\bullet )}(\hat P, \hat P) = 0$ for $i<0$. By the gluing lemma \cite[2.3.3]{Las-Ols05-1} the pair $(P, \sigma )$ is therefore induced by a unique complex on $\X$ which is a perverse sheaf since this can be verified after pulling back to $\mc Y$.  Similarly if $P$ and $P'$ are two perverse sheaves on $\X$ with normalized complexes $\hat P $ and $\hat P'$, then $\ext^i_{\D_c(\X^{\mathbb{N}}, \Lambda _\bullet )}(\hat P, \hat P') = 0$ for $i<0$ and therefore by \cite[2.3.4]{Las-Ols05-1} the functor of morphisms $\hat P\rightarrow \hat P'$ is a sheaf.
\end{proof}

\begin{rem} Using the above argument one can define the category of perverse sheaves on a stack without  defining the $t$--structure.  
\end{rem}

\section{Simple objects}

Let $\X$ be an algebraic stack of finite type over $k$.  Let $\D_c^b(\X, \mathbb{Q}_l)$ denote the category $\D_c^b(\X, \Z_l)\otimes \mathbb{Q}$ (see for example \cite[3.21]{Las-Ols05-2}). The perverse $t$--structure on $\D_c^b(\X, \Z_l)$ defines a $t$--structure on $\D_c^b(\X, \mathbb{Q}_l)$ which we also call the perverse $t$--structure.  An object in the heart of this $t$--structure is called  a perverse $\mathbb{Q}_l$--sheaf.  One check easily that the category of perverse $\mathbb{Q}_l$--sheaves is canonically equivalent to the category $\text{Perv}_{\Z_l}\otimes \Q$, where $\text{Perv}_{\Z_l}$ denotes the category of perverse sheaves of $\Z_l$--modules. In particular, as in ~\ref{perv-stack},the corresponding fibred category is a stack ($\mathbb{Q}_l$-perverse sheaves can be glued).

In what follows we consider only $\mathbb{Q}_l$--coefficients for some $l$ invertible in $k$.

\begin{rem} Verdier duality interchanges $\Dp^{\leq }(\X, \Q_l)$ and $\Dp ^{\geq 0}(\X, \Q_l)$. Indeed this can be verified on a smooth covering of $\X$ and hence follows from the case of schemes.
\end{rem}

\begin{thm}[stack version of {\cite[4.3.1]{BBD82}}]\label{irreduciblethm} \noindent {\rm (i)} In the category of perverse sheaves on $\X$, every object is of finite length.  The category of perverse sheaves is artinian and noetherian.

\noindent {\rm (ii)} Let $j:\mc V\hookrightarrow \X$ be the inclusion of an irreducible substack such that $(\mc V\otimes _k\bar k)_{\text{\rm red}}$ is smooth.  Let $L$ be a smooth $\mathbb{Q}_l$--sheaf on $\mc V$ which is irreducible in the category of smooth $\mathbb{Q}_l$--sheaves on $\mc V$.  Then $j_{!*}(L[\text{\rm dim}(\mc V)])$ is a simple perverse sheaf on $\X$ and every simple perverse sheaf is obtained in this way.
\end{thm}

\begin{proof} Statement (i) can be verified on a quasi--compact smooth covering of $\X$ and hence follows from the case of schemes \cite[4.3.1 (i)]{BBD82}.  

For (ii) note first that if $\X$ is irreducible and smooth, $L$ is a smooth $\Q_l$--sheaf on $\X$, $j:\mc U\hookrightarrow \X$ is a noempty open substack, then the perverse sheaf $F:= L[\text{dim}(\X)]$ satisfies $F = j_{!*}j^*F$.  Indeed it suffices to verify this locally in the smooth topology on $\X$ where it follows from the case of schemes \cite[4.3.2]{BBD82}.

Let $\text{Mod}_{\X}(\Z_l)$ denote the category of smooth adic sheaves of $\Z_l$--modules on $\X$ so that the category of smooth $\Q_l$--sheaves is equal to $\text{Mod}_{\X}(\Z_l)\otimes _{\Z_l}\Q_l$.

\begin{lem} Let $\X$ be a normal algebraic stack of finite type over $k$, and let $j:\mc U\hookrightarrow \mc X$ be a dense open substack.  Then the natural functor
$$
\text{\rm Mod}_{\X}(\Z_l)\rightarrow \text{\rm Mod}_{\mc U}(\Z_l)
$$
is fully faithful and its essential image is closed under subobjects.
\end{lem}
\begin{proof}
Note first that the result is standard in the case when $\X$ is a scheme (in this case when $\X$ is connected the result follows from the surjectivity of the map $\pi _1(\mc U)\rightarrow \pi _1(\X)$). Let $V\rightarrow \X$ be a smooth surjection with $V$ a scheme, and let $U\subset V$ denote the inverse image of $\mc U$.  Also define $V'$ to be $V\times _{\mc X}V$ and let $U'\subset V'$ be the inverse image of $\mc U$.  Assume first that $V'$ is a scheme (in general $V'$ will only be an algebraic space). For any two $F_1, F_2\in \text{Mod}_{\X}(\Z_l)$ we have exact sequences 
$$
0\rightarrow \text{Hom}_{\X}(F_1, F_2)\rightarrow \text{Hom}_{V}(F_1|_V, F_2|_V)\rightrightarrows \text{Hom}_{V'}(F_1|_{V'}, F_2|_{V'})
$$
and
$$
0\rightarrow \text{Hom}_{\mc U}(F_1, F_2)\rightarrow \text{Hom}_{U}(F_1|_V, F_2|_V)\rightrightarrows \text{Hom}_{U'}(F_1|_{U'}, F_2|_{U'}).
$$
From this and the case of schemes the full faithfulness follows.

For the second statement, let $M\in \text{Mod}_{\X}(\Z_l)$ be a sheaf and $L_{\mc U}\subset M|_{\mc U}$ a subobject.  By the case of schemes the pullback $L_U\subset M_V|_U$ to $U$ extends uniquely to a subobject $L_V\subset M_V$.  Moreover, the pullback of $L_V$ to $V'$ via either projection is the unique extension of $L_{U'}$ to a subobject of $M_{V'}$.  It follows that the descent data on $M_V$ induces descent data on $L_V$ restriction to the tautological descent data on $L_U$.  The descended subobject $L\subset M$ is then the desired extension of $L_{\mc U}$.

In all of the above we assumed that $V'$ is a scheme.  This proves in particular the result when $\X$ is an algebraic space.  Repeating the above argument allowing $V'$ to be an algebraic space we then obtain the result for a general stack.
\end{proof}

Tensoring with $\Q_l$ we see that the restriction map
$$
\text{\rm Mod}_{\X}(\Q_l)\rightarrow \text{Mod}_{\mc U}(\Q_l)
$$
is also fully faithful with essential image closed under subobjects.

\begin{lem} Let $\X$ be a normal algebraic stack and $j:\mc U\hookrightarrow \X$ a dense open substack.  If $L$ is a smooth irreducible $\Q_l$--sheaf on $\X$ then the restriction of $L$ to $\mc U$ is also irreducible.
\end{lem}
\begin{proof} Immediate from the preceding lemma.
\end{proof}

\begin{lem} Let $\X$ be a smooth algebraic stack of finite type and $L$ a smooth $\Q_l$--sheaf on $\X$ which is irreducible.  Then the perverse sheaf $F:= L[\text{\rm dim}(\X)]$ is simple.
\end{lem}
\begin{proof}
This follows from the same argument proving \cite[4.3.3]{BBD82} (note that the reference at the end of the proof should be 1.4.25).
\end{proof}

We can now prove \ref{irreduciblethm}.  That the perverse sheaf $j_{!*}F$ is simple follows from \cite[1.4.25]{BBD82} applied to $\mc U\hookrightarrow \overline {\mc U}$ (where $\overline {\mc U}$ is the closure of $\mc U$ in $\mc X$) and \cite[1.4.26]{BBD82} applied to $\overline {\mc U}\hookrightarrow \X$.

To see that every simple perverse sheaf is of this form, let $F$ be a simple perverse sheaf on $\X$.  Then there exists a dense open substack $j:\mc U\hookrightarrow \X$ such that $F_U = L[\text{dim}(\mc U)]$ and such that $(\mc U\otimes _k\bar k)_{\text{red}}$ is smooth over $\bar k$.  By \cite[1.4.25]{BBD82} the map $j_{!*}F_U\rightarrow F$ is a monomorphism whence an isomorphism since $F$ is simple.  This completes the proof of \ref{irreduciblethm}.
\end{proof}

\section{Weights}

In this section we work over a finite field $k=\mathbb{F}_q$.  Fix an algebraic closure $\bar k$ of $k$, and for any integer $n\geq 1$ let $\mathbb{F}_{q^n}$ denote the unique subfield of $\bar k$ with $q^n$ elements.  Following \cite{BBD82} we write objects (e.g. stacks, schemes, sheaves etc.) over $k$ with a subscript $0$ and their base change to $\bar k$ without a subscript.  So for example, $\X_0$ denotes a stack over $k$ and $\X$ denotes $\X_0\otimes _k\bar k$. In what follows we work with $\overline {\Q}_l$--coefficients for some prime $l$ invertible in $k$.

Let $\X_0/k$ be a stack of finite type, and let $\text{Fr}_q:\X\rightarrow \X$ be the Frobenius morphism.  Recall that if $T$ is a $\bar k$--scheme then
$$
\text{Fr}_q:\X(T)=\X_0(T)\rightarrow \X_0(T)=\X(T)
$$
is the pullback functor along the Frobenius morphism of $T$ (which is a $k$--morphism).  We let $\text{Fr}_{q^n}$ denote the $n$--th iterate of $\text{Fr}_q$.  If $x:\Sp (\mathbb{F}_{q^n})\rightarrow \X_0$ is a morphism, we then obtain a commutative diagram
$$
\begin{CD}
\Sp (\bar k)@>\text{Fr}_{q^n}>> \Sp (\bar k)\\
@V\bar xVV @VV\bar xV\\
\X@>\text{Fr}_{q^n}>> \X.
\end{CD}
$$
If $F_0$ is a sheaf on $\X_0$, then the commutativity of this diagram over $\X_0$ defines an automorphism $F_{q^n}^*:F_{\bar x}\rightarrow F_{\bar x}$.

\begin{defn} \noindent {\rm (i)} A sheaf $F_0$ on $\X_0$ is \emph{punctually pure} of weight $w$ ($w\in \Z$) if for every $n\geq 1$ and every $x\in \X_0(\mathbb{F}_{q^n})$ the eigenvalues of the automorphism $F_{q^n}^*:F_{\bar x}\rightarrow F_{\bar x}$ are algebraic numbers all of whose complex conjugates have absolute value $q^{nw/2}$.  

\noindent {\rm (ii)} A sheaf $F_0$ on $\X_0$ is \emph{mixed} if it admits a finite filtration whose successive quotients are punctually pure.  The weights of the graded pieces are called the weights of $F_0$.

\noindent {\rm (iii)} A complex $K\in \D_c^b(\X_0, \overline {\Q}_l)$ is \emph{mixed} if for all $i$ the sheaf $\mc H^i(K)$ is mixed.

\noindent {\rm (iv)} A complex $K\in \D_c^b(\X_0, \overline {\Q}_l)$ is of \emph{weight $\leq w$} if for every $i$ the mixed sheaf $\mc H^i(K)$ has weights $\leq w+i$.

\noindent {\rm (v)} A complex  $K\in \D_c^b(\X_0, \overline {\Q}_l)$ is of \emph{weight $\geq w$} if the Verdier dual of $K$ is of weight $\leq -w$.

\noindent {\rm (vi)} A complex $K\in \D_c^b(\X_0, \overline {\Q}_l)$ is \emph{pure of weight $w$} if it is of weight $\leq w$ and $\geq w$.
\end{defn}

In particular we can talk about a mixed (or pure etc) perverse sheaf.

\begin{thm}[{Stack version of \cite[5.3.5]{BBD82}}] A mixed perverse sheaf $F_0$ on $\X_0$ admits a unique filtration $W$ such that the graded pieces $\text{\rm gr}_i^WF_0$ are pure of weight $i$.  Every morphism of mixed perverse sheaves is strictly compatible with the filtrations.
\end{thm}
\begin{proof} By descent theory (and the uniqueness) it suffices to construct the filtration locally in the smooth topology.  Hence the result follows from the case of schemes.
\end{proof}

The filtration $W$ in the theorem is called the \emph{weight filtration}.

\begin{cor} Any subquotient of a mixed perverse sheaf $F_0$ is mixed. If $F_0$ is mixed of weight $\leq w$ (resp. $\geq w$) then any subquotient is also of weight $\leq w$ (resp. $\geq w$).
\end{cor}
\begin{proof} The weight filtration on $F_0$ induces a filtration on any subquotient whose successive quotients are pointwise pure. This implies the first statement.  The second statement can be verified on a smooth cover of $\X_0$ and hence follows from \cite[5.3.1]{BBD82}.
\end{proof}

One verifies immediately that the subcategory of the category of constructible sheaves on $\X_0$ consisting of mixed sheaves is closed under the formation of subquotients and extensions.  In particular we can define $\D_m^b(\X_0, \overline {\mathbb{Q}}_l)\subset \D_c^b(\X_0, \overline {\Q}_l)$ to be the full subcategory consisting of complexes whose cohomology sheaves are mixed. The category $\D_m^b(\X_0, \overline {\mathbb{Q}}_l)$ is a triangulated subcategory.

\begin{prop} The perverse $t$--structure induces a $t$--structure on $\D_m^b(\X_0, \overline {\mathbb{Q}}_l)$.
\end{prop}
\begin{proof}
It suffices to show that the subcategory $\D_m^b(\X_0, \overline {\mathbb{Q}}_l)\subset \D_c^b(\X_0, \overline {\mathbb{Q}}_l)$ is stable under the perverse truncations $\tau _{\leq 0}$ and $\tau _{\geq 0}$.  This can be verified locally on $\X_0$ and hence follows from the case of schemes.
\end{proof}

%

%
%
%

\end{document}